\documentclass{article}

\def\comment#1{}
\comment{
\oddsidemargin -0.15in %margen(6.3 - 6.5) / 2
\evensidemargin -0.15in %margen (6.3 - 6.5) / 2
\topmargin -0.2in %tope (9.0 - 9.7)
\textwidth=165mm \textheight=228mm
}

\usepackage{amsmath}
\usepackage{color}

\title{Pfaffian Systems of A-Hypergeometric Equations I:
Bases of Twisted Cohomology Groups}
\author{Takayuki Hibi, Kenta Nishiyama, Nobuki Takayama}
\date{June 16, 2014}

\def\QED{ Q.E.D. \par \bigbreak} % use finer Q.E.D. 
\def\CC{ {\bf C} }
\def\ZZ{ {\bf Z} }
\def\RR{ {\bf R} }
\def\pd#1{ \partial_{#1}}
\def\comment#1{ }   % used to comment out a paragraph. 

\newtheorem{corollary}{Corollary}
\newtheorem{theorem}{Theorem}

\newtheorem{remark}{Remark}
\newtheorem{lemma}{Lemma}
\newtheorem{example}{Example}

\begin{document}
\maketitle

{\it Abstract}\/:
We consider bases of Pfaffian systems for $A$-hypergeometric systems.
These are given by Gr\"obner deformations, they also provide bases for twisted cohomology groups.
For a hypergeometric system associated with a class of order polytopes,
these bases have a combinatorial description.
The size of the bases associated with a subclass of the order polytopes
has a growth rate of polynomial order. 

\section{Introduction}

Let $g(x,t) = \sum_{a \in {\cal A}} x_a t^a$,
$ t^a = t_1^{a_1} \cdots t_d^{a_d}$
be a generic sparse polynomial in $t=(t_1, \ldots, t_d)$
with the support on a finite set of points ${\cal A} \subset \ZZ^d$.
The coefficients $x_a$, $a\in {\cal A}$ are denoted by $x_i$, $i=1, \ldots, n$.
The function defined by the integral
$$ \Phi(x) = \int_{C_1} g(x,t)^\alpha t^\gamma dt, \ 
\mbox{ or }\, 
   \Phi(x) = \int_{C_2} \exp(g(x,t)) t^\gamma dt
$$ 
over a cycle $C_i$ in the $t$-space is called
the $A$-{\it hypergeometric function} of $x$
with parameters $\alpha \in \CC$, $\gamma_i \in \CC$
\cite{gkz-euler}, \cite{adolphson}, \cite{esterov-takeuchi}.
It is known that the $A$-hypergeometric function satisfies
a system of linear partial differential equations in $x$,
which is called the $A$-{\it hypergeometric system} or {\it equation}.
The $A$-hypergeometric system is a holonomic system,
and the operators of the system generate a zero-dimensional
ideal in the ring of differential operators with rational function
coefficients (see, e.g., \cite[Chapter 6]{dojo}).
$A$-hypergeometric systems have been studied for the past 25 years
(see, e.g., \cite{gkz}, \cite{gkz-euler}, \cite{SST}), 
and they have applications in various fields.

Let $F$ be a vector-valued function in $x_1, \ldots, x_n$.
We suppose the length of $F$ is $r$ 
and that $F$ is a column vector.
Let $P_i(x)$, $i=1, \ldots, n$, be $r \times r$ matrices satisfying
$$ \frac{\partial P_i}{\partial x_j} + P_i P_j 
  = \frac{\partial P_j}{\partial x_i} + P_j P_i
$$
for all $i \not= j$.
The system of linear differential equations
$$  \frac{\partial F}{\partial x_i} = P_i(x) F, \  i=1, \ldots, n
$$
is called a {\it Pfaffian system}.
We also call the system of linear differential operators
$\frac{\partial}{\partial x_i} - P_i$
a Pfaffian system.
The number $r$ is called the {\it size} or the {\it rank} 
of the Pfaffian system.
For a given zero-dimensional left ideal in the ring of differential operators
with rational coefficients,
it is well known that an associated Pfaffian system can be obtained
by a Gr\"obner basis method (see, e.g., \cite[Appendix]{n3ost2}).
The matrix $P_i$ is an analogy of the companion matrix
for a zero-dimensional ideal in the ring of polynomials.
Some computer algebra systems can perform this translation.
However, in general, this computation is difficult,
and we wish to provide an efficient method for translating
the $A$-hypergeometric system into a Pfaffian system.

Twisted cohomology
groups can be used as a geometric method for finding a Pfaffian system associated with a given
definite integral that contains parameters (see, e.g., the book by Aomoto-Kita \cite{aomoto-kita}).
This approach is as follows: (1) obtain a basis for a twisted
cohomology group, and (2) calculate the Pfaffian system
associated with that basis.
We will use this approach to obtain a Pfaffian system,
and in this paper, we consider step (1).

Gel'fand, Kapranov, and Zelevinsky expressed $A$-hypergeometric functions
with regular singularities as pairings of twisted cycles and twisted cocycles \cite{gkz-euler}.
Esterov and Takeuchi expressed confluent $A$-hypergeometric functions
as pairings of rapid-decay twisted cycles and twisted cocycles \cite{esterov-takeuchi}.
The cohomology groups that come from geometry and are associated with $A$-hypergeometric
systems were discussed by Adolphson and Sperber \cite{adolphson2}.
The next step is to obtain explicit bases for these twisted cohomology groups.
Orlik and Terao provided the $\beta nbc$ bases for the twisted cohomology groups associated with hyperplane
arrangements \cite{orlik-terao}.
Aomoto, Kita, Orlik, and Terao \cite{aomoto-kita-orlik-terao} provided a basis 
for a class of confluent hypergeometric integrals.
In this paper, we will give a computational method for determining the
bases of the twisted cohomology groups associated with generic sparse polynomials
or any $A$-hypergeometric system,
and we will also give a combinatorial method for a class of generic sparse polynomials.

Let $R_n$ be the ring of differential operators with rational
function coefficients in $n$-variables.
The first step in finding a Pfaffian system associated with a zero-dimensional left ideal $I$
in $R_n$ is to obtain a basis for $R_n/I$ as a $\CC(x)$-vector space.
It is well known that a basis can be obtained by computing a Gr\"obner basis
of $I$ in $R_n$.
When the set $\{u_1, \ldots, u_r\}$ is the basis,
there exists a matrix $P_i(x)$ such that
$\partial/\partial x_i  U \equiv P_i U \ {\rm mod}\, I$,
$U = (u_1, \ldots, u_r)^T$.
We can show that $\partial/\partial x_i - P_i$ is a Pfaffian system, and
we call $\{u_1, \ldots, u_r\}$ the {\it basis of the Pfaffian system}.
In Theorem \ref{th:std}, we show that Gr\"obner deformations give bases and 
provide an algorithm that is more efficient than
computing the Gr\"obner basis of $I$ itself.
In Theorem \ref{th:cbasis}, we show that this gives a basis of the twisted cohomology group.

Our theorems are not only useful for computations,
but they also pose interesting theoretical problems in commutative algebra
and combinatorics.
We study the hypergeometric system associated with a class of order polytopes
(see, e.g., \cite{HibiRedBook}).
We prove that the
bases of Pfaffian systems or twisted cohomology groups
have combinatorial descriptions (Theorems \ref{hibi:twochainstandard} and \ref{hibi:bouquettheorem}).
The size of the Pfaffian system associated with a subclass of the order polytopes
has a growth rate of polynomial order (Theorem \ref{hibi:polynomialorder}). 

We will close this introduction by explaining our motivation for this study
from algebraic statistics (see, e.g., \cite{dojo}).
The function $g(x,t)^\alpha t^\gamma/\Phi(x)$
or $\exp(g(x,t))t^\gamma/\Phi(x)$
can be regarded as a probability distribution function on $C_i$
with parameters $x, \alpha, \gamma$ satisfying certain conditions.
This distribution, which we will call the $A$-{\it distribution},
is a generalization of the Beta distribution or the Gamma distribution.
In this context, the function
$\Phi(x)$ is called the normalizing constant of the $A$-distribution.
In \cite{n3ost2}, some new statistical methods were proposed.
These were the holonomic gradient method (HGM) 
and the holonomic gradient descent (HGD).
The HGM is a method for numerically evaluating the normalizing constant,
which is a function of the parameters $x$,
for a given unnormalized probability distribution, 
and the HGD uses the HGM to obtain the maximal likelihood estimate.
The key step for both of these methods is to construct a Pfaffian system
associated with the normalizing constant.
The size of the Pfaffian system determines the complexity of
the HGM and the HGD (see, e.g., \cite{fb2}).
The HGM and HGD lead us to 
the following fundamental goals for applying
$A$-hypergeometric systems to statistics.
\begin{enumerate}
\item Find an efficient method for constructing a Pfaffian system associated with
a given $A$-hypergeometric system.
\item Find a subclass of $A$-hypergeometric systems for which the
associated Pfaffian systems are of moderate size.
\end{enumerate}
We expect that our results will
yield a new class of exponential probability distributions
for which we can efficiently apply the holonomic gradient method (HGM)
and the holonomic gradient descent (HGD).
Construction algorithms for Pfaffian systems, utilizing the results of this paper
and examples of numerical evaluations, will be discussed in the next paper, which is currently
in preparation.

\section{Bases for the Pfaffian System}\label{sec:pf}

We denote by $A=(a_{ij})$ a $d\times n$-matrix whose elements are integers.
We suppose that the set of the column vectors of $A$ spans $\ZZ^d$.
Let $s_1, \ldots, s_d$ be indeterminates.
Following the notation in \cite{SST},
we denote by $H_A(s)$
a left ideal generated  in the Weyl algebra
$$D(s)=\CC(s_1, \ldots, s_d)\langle x_1, \ldots, x_n, \pd{1}, \ldots, \pd{n} \rangle, \ \pd{i}=\partial/\partial x_i
$$
by
\begin{eqnarray}
  &&\sum_{j=1}^n a_{ij}  x_j \partial_{j} - s_i, 
   \qquad(i = 1, \ldots, d)  \label{eq:euler} \\
  &&\prod_{i=1}^n \partial_{i}^{u_{i}} 
  - \prod_{j=1}^n \partial_{j}^{v_{j}} \label{eq:toric} \\
  && \quad\quad (\mbox{with } u, v \in {\bf N}_0^{n} \mbox{ running over all $u, v$ such that } A u = Av ).  \nonumber
\end{eqnarray}
Here, ${\bf N}_0 = \{ 0, 1, 2, \ldots \}$.
We call the ideal generated in $\CC[\pd{1}, \ldots, \pd{n}]$
by the elements of the form (\ref{eq:toric})
the affine toric ideal and denote it by $I_A$.
We denote by $E_i-s_i$ the operator (\ref{eq:euler}).
For complex parameters $\beta_i$, 
the system of linear differential equations
$(E_i - \beta_i) f = 0$\,  ($i=1, \ldots, d$),
$(\pd{}^u - \pd{}^v) f = 0$\, ($Au=Av$)
is called the $A$-hypergeometric system of differential equations
or just the $A$-hypergeometric system.
We will sometimes call the ideal $H_A(s)$ 
the $A$-hypergeometric system (with indefinite parameters).

Let $R_n$ be the ring of differential operators
with rational function coefficients
\begin{equation} \label{eq:defOfR}
  \CC(s,x) \langle \pd{1}, \ldots, \pd{n} \rangle.
\end{equation}
We are interested in bases of $R_n/(R_n H_A(s))$
as the vector space over the field $\CC(s,x)$.
Any basis of the vector space yields an associated Pfaffian system
or an integrable connection associated with $H_A(s)$.
Let $u_1, \ldots, u_r$ be a basis of $R_n/(R_n H_A(s))$.
For $u_j$, there exist rational functions $p_{ij}^k \in \CC(s,x)$ 
such that 
$\pd{i} u_j \equiv \sum_{k=1}^r p_{ij}^k u_k \, {\rm mod} \, R_n H_A(s)$.
The action of a differential operator $u$ to a function $F$ is denoted by
$u \bullet F$.
The system of differential equations
$ \pd{i}\bullet F = (p_{ij}^k \,|\, 1 \leq j,k \leq r) F$,
where $F$ is a vector valued function of size $r$,
is called a Pfaffian system,
and $\{ u_i \}$  is called a basis of the Pfaffian system.

Bases can be described by those of simpler quotients,
for which the denominator ideals are nothing but Gr\"obner deformations of $H_A(s)$,
as in the following theorem.

\begin{theorem} \label{th:std}
Let $w \in \ZZ^n$ be a generic weight vector for the affine toric ideal $I_A$
such that ${\rm deg}\, {\rm in}_w(I_A) = {\rm deg}\, I_A$.
Let $u_1, \ldots, u_r$ be a monomial basis of $R_n/(R_n J)$,
where the left ideal $J$ is generated by
${\rm in}_w(I_A)$ and
$E_i-s_i$, $i=1, \ldots, d$ in $R_n$.
Then, the set $\{ u_1, \ldots, u_r \}$ is a basis of the vector space
$R_n/(R_n H_A(s))$.
\end{theorem}

{\it Proof}\/.
We denote by $r$ the normalized volume of $A$.
Since the $s_i$ are indeterminate, the holonomic rank of $J$ and $H_A(s)$
are $r$ by Adolphson's theorem (see, e.g., \cite{adolphson}, \cite{SST}).
In other words, we have
${\rm dim}_{\CC(s,x)} R_n/(R_n H_A(s)) = r$
and
${\rm dim}_{\CC(s,x)} R_n/J = r$.

We may assume that the $u_i$ are expressed as monomials in terms of 
Euler operators $\theta_j = x_j \pd{j}$.
When we regard $J$ as a system of linear differential equations,
it has $r$ linearly independent solutions of the form $x^\rho$,
where $\rho \in \CC(s)^n$.
We denote them by $g_i=x^{\rho(i)}$, $i=1, \ldots, r$.
Since the $g_i$ are linearly independent solutions, the Wronskian
determinant ${\rm det}(u_i \bullet g_j)$ is not identically 
equal to $0$.
The solution $g_j$ can be extended to a solution $f_j$ of $H_A(s)$
such that $g_j$ is the leading monomial of $f_j$ with respect to the
weight vector $w$
(see, e.g., \cite[Chapters 2 and 3]{SST}).
The series $f_j$ is expressed as
$
 f_j = g_j \sum_{\ell \in M_j} C_\ell x^\ell
$, $C_0 = 1$,
where $M_j$ denotes the set of lattice points in a cone and
$C_\ell$ is a constant belonging to $\CC(s)$.
The series converges  in the space of convergent power series 
$g_j \cdot {\cal O}(U)\{M_j\}$, 
where $U$ is an open set in the $s$-space
and ${\cal O}(U)$ is the space of holomorphic functions on $U$
\cite{ohara-takayama}.
We replace $x_i$ by $x_i t^{w_i}$ for all $i$ in $f_j$ and
denote by $x t^w$ the vector $(x_1 t^{w_1}, \ldots, x_n t^{w_n})$.
From the construction algorithm of $f_j$, we may assume that
$f_j(xt^w) = g_j(xt^w) ( 1 + O(t))$ 
when $t \rightarrow 0$ as a function of $t$ 
when $x$ is fixed and $s$ lies in $U$.

Let us prove $W={\rm det}(u_i \bullet f_j) \not\equiv 0$.
We denote by $u_i(\rho(j))$ 
the constant $ (u_i \bullet x^{\rho(j)}) /x^{\rho(j)}$.
Under this notation, we have
$x^{-\rho(j)} u_i \bullet g_j = u_i(\rho(j)) $ and
\begin{equation} \label{eq:series1}
  (x^{-\rho(j)} (u_i \bullet f_j))(x t^w)
 = \sum_{\ell \in M_j} u_i(\rho(j)+\ell) C_\ell x^\ell t^{\ell w}.
\end{equation}
Note that $\ell w > 0$ for $\ell \not=0$ and $\ell \in M_j$.
Therefore, we have
\begin{equation}
  {\rm det}(x^{-\rho(j)} u_i \bullet f_j)(x t^w)
= {\rm det}(x^{-\rho(j)} u_i \bullet g_j)(x t^w) + O(t)
\end{equation}
from (\ref{eq:series1}) when
$x$ is fixed and $t \rightarrow 0$.
This implies that the Wronskian determinant 
${\rm det}(u_i \bullet f_j)
= \left(\prod_j x^{\rho(j)}\right) {\rm det}(x^{-\rho(j)} u_i \bullet f_j)$
is not identically equal to $0$.
Therefore the $u_i$ are linearly independent in $R_n/(R_n H_A(s))$.
\QED
%%Ref: 	@s/2012/10/11-my-note-basis-of-pfaffian.pdf  std mon of in I_A + E gives a std mon of H_A.

Let $M$ be a monomial ideal in $\CC[\pd{}]$.
When $M$ is generated by $\pd{}^\alpha$,
the distraction ${\widetilde M} \subset {\bf C}[\theta]$ is generated by
$\prod_{i=1}^n \theta_i (\theta_i-1) \cdots (\theta_i-\alpha_i+1)$,
where $\theta_i = x_i \pd{i}$ \cite[p.68]{SST}.
Let $M={\rm in}_w(I_A)$.
Then, the ideal $J$ in Theorem \ref{th:std} is generated by 
${\widetilde M}$ and $\sum_{j=1}^n a_{ij} \theta_j - s_i$, $i=1, \ldots, d$
\cite[Sec. 2.3, Prop. 3.1.5]{SST}.
This leads us to the following corollary.
\begin{corollary} \label{cor:std}
Retain the assumptions of Theorem \ref{th:std}.
The set of the monomial basis of $\CC(s)[\theta]/{\tilde J}$,
where
${\tilde J}$ is the ideal generated by
${\widetilde M}$ and 
$\sum a_{ij} \theta_j-s_i$, $i=1, \ldots, d$ in the polynomial ring
$\CC(s)[\theta]$,
gives a basis of $R_n/R_n H_A(s)$
by the replacement $\theta_i = x_i \pd{i}$. 
\end{corollary}

\section{Bases of Twisted Cohomology Groups} \label{sec:twisted}

Let $A_1 = (a_1, \ldots, a_{n_1}), \ldots,
 A_k = (a_{n_{k-1}+1}, \ldots, a_{n_k})$, $a_i \in \ZZ^m$.
To each matrix $A_j$, we associate a generic sparse polynomial in $t$
\begin{equation} \label{eq:fdef}
 f_j(x,t) = \sum_{i=n_{j-1}+1}^{n_j} x_i t^{a_i},
\end{equation}
where $t^b = \prod_{i=1}^m t_i^{b_i}$.
For parameters $\alpha_1, \ldots, \alpha_k$ and
$\gamma_1, \ldots, \gamma_m$,
we consider the integral
\begin{equation} \label{eq:kernel}
\Phi(\alpha,\gamma; x)=\int_C P(x,t) dt_1 \cdots dt_m, \ 
P(x,t)=\prod_{j=1}^k f_j(x,t)^{\alpha_j} t^\gamma 
\end{equation}
for a suitable twisted cycle $C$.
The function $\Phi$ is satisfied by the $A$-hypergeometric system
for 
\begin{equation} \label{eq:aaa}
A = \left( \begin{array}{cccccccccccc}
             1 & \cdots & 1  & 0 & \cdots & 0   & & & &    0 & \cdots & 0 \\ 
             0 & \cdots & 0  & 1 & \cdots & 1   & & & &    0 & \cdots & 0 \\ 
             0 & \cdots & 0  & 0 & \cdots & 0   & & & &    0 & \cdots & 0 \\ 
          \cdot& \cdots &\cdot&\cdot&\cdots&\cdot&& & & \cdot& \cdots &\cdot\\ 
          \cdot& \cdots &\cdot&\cdot&\cdots&\cdot&\cdot&\cdot&\cdot & \cdot& \cdots &\cdot\\ 
          \cdot& \cdots &\cdot&\cdot&\cdots&\cdot&& & & \cdot& \cdots &\cdot\\ 
             0 & \cdots & 0  & 0 & \cdots & 0   & & & &    1 & \cdots & 1 \\ 
            a_1& \cdots&a_{n_1}&a_{n_1+1}&\cdots&a_{n_2}& & & &a_{n_{k-1}+1}& \cdots &a_{n_k}\\ 
\end{array} \right)
\end{equation}
and 
$\beta = ( \alpha_1, \ldots, \alpha_k, -\gamma_1-1, \ldots, -\gamma_m-1)^T$,
where we assume that the rank of $A$ is maximal
\cite{gkz-euler}.
Set $P' = P {|_{\alpha=\gamma=1}}$.
Let $n = \sum_{i=1}^k n_i$, and
define a projection $p$ by 
$$ p\ :\  {\bf C}^{m+n} \setminus V(P') \ni (x,t) \mapsto x \in {\bf C}^n. $$
We regard $P$ as a function in $t=(t_1, \ldots, t_{m})$ with the parameter vector $x$.
Define the connection $\nabla$ with rational function coefficients by
\begin{equation} \label{eq:conn}
\nabla = d + \sum_{j=1}^{m} \left( \frac{\partial P}{\partial t_j}/P \right) dt_j
\end{equation}
where $d$ is the exterior derivative with respect to the variables $t_1, \ldots, t_{m}$.

\begin{theorem}  \label{th:cbasis}
Assume that the matrix $A$ is expressed as (\ref{eq:aaa}).
Let $\alpha, \gamma$ be generic parameters,
and let $\{u_1, \ldots, u_r\}$ be a basis as given in Theorem \ref{th:std}.
Then the set of rational expressions
$\{ (u_i \bullet P)/P \}$ is a basis of the twisted cohomology group
$H^m (p^{-1}(x),\nabla)$ when $x$ lies outside of an analytic set.
\end{theorem}

{\it Proof}\/.
It follows from the local triviality theorem \cite[5.1 Corollaire]{verdier}
that the projection $p$
is a locally trivial map on a Zariski open subset $U$ of ${\bf C}^n$.
Take a point $x_0$ in $U$.
Then the inverse image $p^{-1}(U')$ of a small neighborhood $U'$ of $x_0$
is isomorphic, as a smooth manifold, to the direct product of $p^{-1}(x_0) \times U'$.
Therefore, we can form a basis from the twisted homology group 
$ H_m (p^{-1}(x), {\cal P}_x)$, $x \in U'$
of the form
$ \sum c_i \Delta_i \otimes P$,
where $c_i$ is a constant that does not depend on $x$ and $\Delta_i$ is a smooth simplex
that does not depend on $x$.
Here, ${\cal P}_x$ is the local system defined by $P$ at $x$. Note that
${\cal P}_x$ and ${\cal P}_{x'}$ are isomorphic for any $x, x' \in U'$.
%%Note: therefore, we can take $c_i$ and $\Delta_i$ which do not depend on x. @s/2013/10/10-my-note-*.pdf

Let $\{ u_1, \ldots, u_r \}$ be a basis of $R_n/(R_n H_A(s))$ as given in Theorem \ref{th:std}.
For generic parameters $\alpha$ and $\gamma$ and a twisted cycle
$ C_j = \sum_k c_{jk} \Delta_{jk} \otimes P$, where $c_{jk}$ and $\Delta_{jk}$
does not depend on the parameter $x$,
the integral
$$ u_i \Phi(\alpha,\gamma;x) = \sum c_{jk} \int_{\Delta_{jk}} u_i \bullet P dt
  = \sum c_{jk} \int_{\Delta_{jk}} \frac{u_i \bullet P}{P} P dt
$$
can be regarded as a pairing $\langle \varphi_i,C_j\rangle$ of the twisted cocycle 
$ \varphi_i = \frac{u_i \bullet P}{P}dt \in H^m (p^{-1}(x),\nabla)$ and
the twisted cycle $C_j$.
Since the matrix-valued function $(\langle \varphi_i, C_j \rangle)$
is a fundamental set of solutions of the Pfaffian system for the $A$-hypergeometric 
system $H_A(\beta)$,
its determinant does not vanish out of an analytic set.
This implies that the pairing of
$H^m (p^{-1}(x),\nabla) \times H_m(p^{-1}(x),{\cal P}_x)$ 
is perfect out of the analytic set in the $x$ space.
Thus, the set $u_i \bullet P/P $ is a basis of the twisted cohomology group.
\QED

%%Todo We note that a generalization of the Theorem \ref{th:twisted} hold  --> done 2013.11.21.
%% for more general $A$ under the condition ???  Orlik-Terao's results.
%%Ref: 2012/12/22-my-note-mhg-conf-f1-gkz-advance.pdf, *takeuchi*, H_c and H, vanishing.

\begin{remark} \label{rem:rapid_decay} \rm
Our Theorem \ref{th:cbasis} can be generalized to a matrix $A$ for which the toric ideal $I_A$ is not necessarily 
a homogeneous ideal.
It follows from Esterov and Takeuchi \cite{esterov-takeuchi} that we can form a basis from the rapid-decay homology cycles
for which the support does not depend locally on the parameter $x$.
Therefore, we can make an analogous argument for the case of twisted cycles, and
the perfectness theorem by Hien \cite{hien} proves Theorem \ref{th:cbasis}
for any $A$ that defines an $A$-hypergeometric system.
\end{remark}

\begin{figure}[tb]
\setlength{\unitlength}{0.3mm}
\begin{center}
\begin{picture}(100,105)(0,15)
\put(50,25){\circle*{2}}
\put(50,20){1}
\put(25,50){\circle*{2}}
\put(25,45){2}
\put(50,50){\circle*{2}}
\put(50,45){3}
\put(75,50){\circle*{2}}
\put(75,45){4}
\put(25,75){\circle*{2}}
\put(25,70){5}
\put(50,75){\circle*{2}}
\put(50,70){6}
\put(75,75){\circle*{2}}
\put(75,70){7}
\put(50,100){\circle*{2}}
\put(50,95){8}

\put(50,25){\line(-1,1){25}}
\put(50,25){\line(0,1){25}}
\put(50,25){\line(1,1){25}}

\put(25,50){\line(0,1){25}}
\put(25,50){\line(1,1){25}}

\put(50,50){\line(-1,1){25}}
\put(50,50){\line(1,1){25}}

\put(75,50){\line(-1,1){25}}
\put(75,50){\line(0,1){25}}
%\put(75,50}{\line(1,1){25}}

\put(50,100){\line(-1,-1){25}}
\put(50,100){\line(0,-1){25}}
\put(50,100){\line(1,-1){25}}
\end{picture}
\end{center}
\caption{$C_{111}$}  \label{fig:c1}
\end{figure}

\begin{example} \rm \label{ex:c11}
Consider the matrix
$$
A=
\bordermatrix{
  & x_1 & x_2 & x_3 & x_4 & x_5 & x_6 & x_7 & x_8 \cr
  & 1 & 1 & 1 & 1 & 1 & 1 & 1 & 1 \cr
  & 0 & 1 & 0 & 0 & 1 & 1 & 0 & 1 \cr
  & 0 & 0 & 1 & 0 & 1 & 0 & 1 & 1 \cr
  & 0 & 0 & 0 & 1 & 0 & 1 & 1 & 1 \cr
}
$$
(which are the vertices of the order polytope associated to the distributive lattice of Figure \ref{fig:c1}; see Section \ref{sec:orderpoly}).
The basis given by Theorem \ref{th:std} is
$\{
 1, \pd{5}, \pd{6}, \pd{7}, \pd{8}, \pd{8}^2
\}
$.
It is determined by computing a Gr\"obner basis for the ideal $J$ or ${\tilde J}$
by a computer.
Note that this computation is easier than computing the Gr\"obner basis
of $H_A(\beta)$.
The corresponding basis of the twisted cohomology group is
$\{
 1, \frac{t_1 t_2 dt}{Q} , \frac{t_1 t_3dt}{Q}, \frac{t_2 t_3dt}{Q},
 \frac{t_1 t_2 t_3dt}{Q}, \frac{(t_1 t_2 t_3)^2dt}{Q^2}
\},
$
where $Q=x_1+x_2 t_1 + x_3 t_2 + x_4 t_3 + x_5 t_1 t_2 + x_6 t_1 t_3 + x_7 t_2 t_3 + x_8 t_1 t_2 t_3$
and $dt = dt_1 dt_2 dt_3$.
\end{example}

\begin{example}\rm \label{ex:confc11}
Let
$$
A'=
\bordermatrix{
   & x_2 & x_3 & x_4 & x_5 & x_6 & x_7 & x_8 \cr
   & 1 & 0 & 0 & 1 & 1 & 0 & 1 \cr
   & 0 & 1 & 0 & 1 & 0 & 1 & 1 \cr
   & 0 & 0 & 1 & 0 & 1 & 1 & 1 \cr
}.
$$
The hypergeometric system associated with $A'$ is the confluent system
of the previous example.
Put $Q=x_2 t_1 + x_3 t_2 + x_4 t_3 + x_5 t_1 t_2 + x_6 t_1 t_3 + x_7 t_2 t_3 + x_8 t_1 t_2 t_3$
and $dt = dt_1 dt_2 dt_3$.
Then the integral $\int_C \exp(Q) t^\gamma dt$, where $C$ is a rapid decay
cycle, is a solution.
The toric ideal $I_{A'}$ is obtained formally
by setting $\pd{1} = 1$ in $I_A$ of Example \ref{ex:c11}.
%%Prog/cbase.rr  ctest5();
A basis given by Theorem \ref{th:std} is
$\{
 1, \pd{5}, \pd{6}, \pd{7}, \pd{8}, \pd{8}^2
\}
$.
The corresponding basis of the twisted cohomology group is
$\{
 1, t_1 t_2 dt , t_1 t_3dt, t_2 t_3dt,
 t_1 t_2 t_3dt, (t_1 t_2 t_3)^2dt
\}
$.
\end{example}

These two examples illustrate that the results in Sections \ref{sec:pf} and \ref{sec:twisted} give 
a general method
for determining the bases of twisted cohomology groups by computing a set of standard monomials
with a Gr\"obner basis of $J$ (Theorem \ref{th:std}) or ${\tilde J}$ (Corollary \ref{cor:std}).

\section{$A$-hypergeometric Systems for Order Polytopes}\label{sec:orderpoly}
% Let $P$ be a finite partially ordered set (\cite[p.\,000]{StanleyEC}) 
% of $|p| = d$ and $L = {\mathcal J}(P)$ the distributive lattice (\cite[p.\,000]{StanleyEC}) 
% consisting of all poset ideals (\cite[p.\,000]{StanleyEC}) of $P$ ordered by inclusion.
% The order polytope ${\mathcal O}(P)$ (\cite[p.\,000]{StanleyEC}) arising from $P$ 
% has been studied from viewpoints of both commutative algebra and combinatorics.  
% Especially, in \cite{Hibi}, the toric ideal of ${\mathcal O}(P)$ 
% and its Gr\"obner basis is studied. 
For some classes of generic sparse polynomials or $A$, we can calculate by hand the set of standard monomials
for $J$ or $\tilde J$.

First, recall the order polytope % ${\mathcal O}(P)$  
of a finite partially ordered set (\cite[p.\,115]{HibiRedBook}).
Let $P = \{a_{1}, \ldots, a_{d}\}$ be a finite partially ordered set
with $|P| = n$. 
A {\em poset ideal} of $P$ is a subset $\alpha$ of $P$ such that 
if $a \in \alpha$, $b \in P$, and $b \leq a$, then $b \in \alpha$.
Thus in particular the empty set and $P$ itself are poset ideals of $P$.
Let ${\mathcal J}(P)$ denote the distributive lattice (\cite[p.\,118]{HibiRedBook}) 
consisting of all poset ideals of $P$, ordered by inclusion.
For example, if $P$ is the disjoint union of two chains of length $2$ and length $3$
shown in Figure \ref{fig:l23}, then $L = {\mathcal J}(P)$ is the distributive lattice
shown in Figure \ref{fig:s34}.

\begin{figure}[tbp]
\begin{minipage}{0.49\hsize}
\setlength{\unitlength}{0.3mm}
\begin{center}
\begin{picture}(100,100)(0,15)
\put(33,10){\line(0,1){50}}
\put(66,10){\line(0,1){75}}
\multiput(33,10)(0,25){3}{\circle*{3}}
\multiput(66,10)(0,25){4}{\circle*{3}}
\end{picture}
\end{center}
\caption{A poset $P$}  \label{fig:l23}
\end{minipage}
%%%%
\begin{minipage}{0.49\hsize}
\setlength{\unitlength}{0.3mm}
\begin{center}
\begin{picture}(100,100)(0,0)
\multiput(45,0)(15,15){5}{\line(-1,1){45}}
\multiput(45,0)(-15,15){4}{\line(1,1){60}}
\multiput(45,0)(15,15){5}{\multiput(0,0)(-15,15){4}{\circle*{3}}}
\end{picture}
\end{center}
\caption{Distributive lattice of ${\mathcal J}(P)$}  \label{fig:s34}
\end{minipage}
\end{figure}

Let ${\bf e}_{1}, \ldots, {\bf e}_{d}$ denote the standard unit coordinate vectors of
${\bf R}^{d}$.  If $\beta$ is a subset of $P$, then we write $w_{\beta}$ for
the $(0,1)$-vector $\sum_{a_{i}\in \beta} {\bf e}_{i} \in {\bf R}^{d}$.
The {\em order polytope} ${\mathcal O}(P) \subset {\bf R}^{d}$ of $P$ 
is the convex hull of the finite set 
$\{ w_{\alpha} \, : \, \alpha \in {\mathcal J}(P) \}$.
Its dimension is $\dim {\mathcal O}(P) = d$.   

Let $K = {\bf C}(\{\xi_{\alpha}\}_{\alpha \in {\mathcal J}(P)})$ 
denote the rational function field
in $|{\mathcal J}(P)|$ variables over % the complex number field 
${\bf C}$.
Let $A = K[t_{1}, \ldots, t_{d},s]$ denote the polynomial ring in $n$ variables 
over $K$. If $\beta$ is a subset of $P$, then we write $u_{\beta}$ for the
square-free  monomial $\prod_{a_{i} \in \beta} t_{i}s$.
Let $K[{\mathcal O}(P)]$ denote the subalgebra of $A$ that is generated by
those square-free  monomials $u_{\beta}$ with $\beta \in {\mathcal J}(P)$.
The semigroup ring $K[{\mathcal O}(P)]$ was introduced in \cite{Hibi}.
We call $K[{\mathcal O}(P)]$ the {\em toric ring} of ${\mathcal O}(P)$. 
The Krull dimension of $K[{\mathcal O}(P)]$ is $d + 1$.

Let $K[\{x_{\alpha}\}_{\alpha \in {\mathcal J}(P)}]$ denote the polynomial ring
in $|{\mathcal J}(P)|$ variables over $K$, and define the surjective ring homomorphism
$\pi : K[\{x_{\alpha}\}_{\alpha \in {\mathcal J}(P)}] \rightarrow K[{\mathcal O}(P)]$
by setting $\pi(x_{\alpha}) = u_{\alpha}$.  Its kernel $I_{{\mathcal O}(P)}$
is called the {\em toric ideal} of ${\mathcal O}(P)$.  
It is known \cite{Hibi} that $I_{{\mathcal O}(P)}$ is generated by quadratic binomials
\begin{eqnarray}
\label{hibirelation}
x_{\alpha}x_{\beta} - x_{\alpha \wedge \beta}x_{\alpha \vee \beta},
\end{eqnarray}
such that $\alpha$ and $\beta$ are incomparable in 
the distributive lattice ${\mathcal J}(P)$.
We fix an ordering $<$ of variables of $K[\{x_{\alpha}\}_{\alpha \in {\mathcal J}(P)}]$
with the property that
if $\alpha > \beta$ in ${\mathcal J}(P)$, then $x_{\alpha} < x_{\beta}$.
Let $<_{\rm rev}$ be the reverse lexicographic order on  
$K[\{x_{\alpha}\}_{\alpha \in {\mathcal J}(P)}]$
induced by the ordering $<$.  In \cite{Hibi} it was shown that 
the set of binomials (\ref{hibirelation}) is the reduced Gr\"obner basis of
$I_{{\mathcal O}(P)}$ with respect to $<_{\rm rev}$.
Thus ${\rm in}_{<_{\rm rev}}(I_{{\mathcal O}(P)})$ is generated by those square-free  
quadratic monomial $x_{\alpha}x_{\beta}$ such that 
$\alpha$ and $\beta$ are incomparable in ${\mathcal J}(P)$.

Let
\[
\theta_{i} = \sum_{a_{i}\in \alpha} \xi_{\alpha}x_{\alpha} - \eta_{i}, \, \, \, \, \, \, \, \, \, \, 
1 \leq i \leq d,
\]
and let
\[
\theta_{0} = \sum_{\alpha \in {\mathcal J}(P)} \xi_{\alpha}x_{\alpha} - \eta_{0},
\]
where $\eta_i \in K$.
It then follows that the sequence $(\theta_{0}, \theta_{1}, \ldots, \theta_{d})$ 
is a system of parameters of both the residue rings 
$K[\{x_{\alpha}\}_{\alpha \in {\mathcal J}(P)}]/I_{{\mathcal O}(P)}$ 
and 
$K[\{x_{\alpha}\}_{\alpha \in {\mathcal J}(P)}]/{\rm in}_{<_{\rm rev}}(I_{{\mathcal O}(P)})$.
The fundamental goal is to find
a $K$-basis of the zero-dimensional 
residue ring 
\begin{eqnarray}
\label{hibi:0-dimensional}
K[\{x_{\alpha}\}_{\alpha \in {\mathcal J}(P)}]/({\rm in}_{<_{\rm rev}}(I_{{\mathcal O}(P)}),
\theta_{0}, \theta_{1}, \ldots, \theta_{d}).
\end{eqnarray}
%By virtue of Theorem \ref{th:std_commutative}, the set of standard monomials of (\ref{hibi:0-dimensional})
%is a $K$-basis of
%\[
%K[\{x_{\alpha}\}_{\alpha \in {\mathcal J}(P)}]/(I_{{\mathcal O}(P)},
%\theta_{0}, \theta_{1}, \ldots, \theta_{d}).
%\] 
In general, however, this is difficult. 
When $P$ can be decomposed into two chains, a complete answer can be found, as shown below.
For example, the poset $P$ of Figure \ref{fig:l23m} can be decomposed into
the chains $a_{1} < a_{2} < a_{3}$ and $b_{1} < b_{2} < b_{3} < b_{4}$.   

\begin{figure}[tbp]
\begin{minipage}{0.49\hsize}
\setlength{\unitlength}{0.7mm}
\begin{center}
\begin{picture}(30,35)(0,-5)
\put(10,0){\line(0,1){20}}
\put(20,0){\line(0,1){30}}
\multiput(10,0)(0,10){3}{\circle*{3}}
\multiput(20,0)(0,10){4}{\circle*{3}}
\put(0,0){$a_1$}
\put(0,10){$a_2$}
\put(0,20){$a_3$}
\put(25,0){$b_1$}
\put(25,10){$b_2$}
\put(25,20){$b_3$}
\put(25,30){$b_4$}
\put(10,10){\line(1,2){10}}
\put(10,20){\line(1,-2){10}}
\end{picture}
\end{center}
\caption{$P$}  \label{fig:l23m}
\end{minipage}
%%%%
\begin{minipage}{0.49\hsize}
\setlength{\unitlength}{0.3mm}
\begin{center}
\begin{picture}(100,100)(0,0)
\multiput(45,0)(15,15){4}{\line(-1,1){30}}
\multiput(45,0)(-15,15){3}{\line(1,1){45}}
\multiput(45,0)(15,15){4}{\multiput(0,0)(-15,15){3}{\circle*{3}}}
\multiput(15,60)(15,15){4}{\circle*{3}}
\multiput(15,60)(15,15){3}{\line(1,1){15}}
\multiput(15,60)(15,15){4}{\line(1,-1){15}}
\put(75,90){\circle*{3}}
\put(75,90){\line(-1,-1){15}}
\end{picture}
\end{center}
\caption{${\mathcal J}(P)$}  \label{fig:s34m}
\end{minipage}
\end{figure}

Now, suppose that a finite poset $P$ can be decomposed into
two chains $C_{p} : a_{1} < \cdots < a_{p}$
of length $p-1$ 
and $C_{q} : b_{1} < \cdots < b_{q}$ of length $q-1$,
where $p \geq 1$ and $q \geq 1$.
Let ${\mathcal A}$ denote the set of those pairs $(i,j)$,
where $0\leq i \leq p$ and $0\leq j \leq q$, for which
$\{a_{1},\ldots,a_{i}, b_{1},\ldots,b_{j}\}$ is a poset ideal of $P$.
In particular, $(0,0), (p, q) \in {\mathcal A}$.   
When $(i,j) \in {\mathcal A}$, 
we write $\alpha_{i,j}$ for $\{a_{1},\ldots,a_{i}, b_{1},\ldots,b_{j}\}$.
For example, $\alpha_{0,0} = \emptyset$ and $\alpha_{p,q} = P$. 
We then have $L = \{ \, \alpha_{i,j} \, : \, (i,j) \in {\mathcal A} \, \}$.
When $(i,j) \in {\mathcal A}$, we write
$\xi_{i,j}$ for $\xi_{\alpha_{i,j}}$ and
$x_{i,j}$ for $x_{\alpha_{i,j}}$.
% Thus we are working with the polynomial ring 
% then we may identify $\alpha_{i,j} \in {\mathcal J}(P)$ with the point $(i,j) \in {\bf R}^{2}$.
% \[
% K[\{x_{\alpha}\}_{\alpha \in {\mathcal J}(P_{p,q})}]=
% K[\{x_{\alpha_{i,j}}\}_{1\leq i \leq p, 1\leq j \leq q}]=
% K[\{x_{i,j}\}_{1\leq i \leq p, 1\leq j \leq q}],
% \]
Let
\[
\theta_{i*} = \sum_{i \leq k \leq p, \, 0\leq j \leq q, \, (k,j) \in {\mathcal A}} 
\xi_{k,j}x_{k,j} - \eta_{i*}, \, \, \, \, \, \, \, \, \, \, 
0 \leq i \leq p
\]
and
\[
\theta_{*j} = \sum_{0 \leq i \leq p, \, j\leq \ell \leq q, \, (i,\ell) \in {\mathcal A}} 
\xi_{i,\ell}x_{i,\ell} - \eta_{*j}, \, \, \, \, \, \, \, \, \, \, 
0 \leq j \leq q.
\]
In particular,
\[
\theta_{0*} = \theta_{*0} = \sum_{(i,j) \in {\mathcal A}} 
\xi_{i,j}x_{i,j} - \eta_{0},
\]
with $\eta_{0} = \eta_{0*} = \eta_{*0}$.
Let $K[{\bf x}]=K[\{x_{i,j}\}_{0\leq i \leq p, \, 0\leq j \leq q, \, 
(i,j) \in {\mathcal A}}]$ and
\begin{equation} \label{eq:J}
J = ( \, {\rm in}_{<_{\rm rev}}(I_{{\mathcal O}(P)}), \,
\{\theta_{i*}\}_{0\leq i \leq p}, \, \{\theta_{*j}\}_{0\leq j \leq q} \, ),
\end{equation}
where
\[
{\rm in}_{<_{\rm rev}}(I_{{\mathcal O}(P)}) = 
(\{ \, x_{i,j}x_{k,\ell} \, : \, i < k, \, \ell < j, \, 
(i,j) \in {\mathcal A}, \, (k, \ell) \in {\mathcal A} \, \}).
\]
Then the residue ring (\ref{hibi:0-dimensional}) is
$K[{\bf x}]/J$.
Let $<_{\rm rev}$ denote the reverse lexicographic order on $K[{\bf x}]$
induced by the ordering of the variables, as follows:
$
x_{i,j}>x_{k,\ell}
$
if either $i + j < k + \ell$ or $i + j = k + \ell$ with $i > k$.

\begin{lemma}
\label{hibi:old}
In $K[{\bf x}]/{\rm in}_{<_{\rm rev}}(J)$, 
\[
x_{i,j} x_{i,j'} \, = \, x_{i,j} x_{i',j} \, = \, 0,
\]
% together with $x_{i,0}=x_{0,j}=0$,
where $(i,j), (i,j')$ and $(i',j)$ belong to ${\mathcal A}$.
\end{lemma}

{\it Proof}\/.  
Let $i < i'$. Then
% \[
% \theta_{i'*}x_{i,j} - x_{i,j} (\sum_{i' \leq k, \, j \leq \ell, \, 
% (k, \ell) \in {\mathcal A}}
% \xi_{k,\ell}x_{k,\ell})
% \]
\[
\theta_{i'*}x_{i,j} - x_{i,j} \Big(\Big(\sum_{i' \leq k, \, j \leq \ell, \, 
(k, \ell) \in {\mathcal A}}
\xi_{k,\ell}x_{k,\ell} \Big) + b_{i'*}\Big)
\] 
belongs to ${\rm in}_{<_{\rm rev}}(I_{{\mathcal O}(P)})$.  Hence
% \[
% x_{i,j} (\sum_{i' \leq k, \, j \leq \ell, \, (k, \ell) \in {\mathcal A}}
% \xi_{k,\ell}x_{k,\ell})
% \]
\[
x_{i,j} \Big(\Big(\sum_{i' \leq k, \, j \leq \ell, \, 
(k, \ell) \in {\mathcal A}}
\xi_{k,\ell}x_{k,\ell} \Big) + b_{i'*}\Big)
\]
belongs to $J$.  
Thus its initial monomial $x_{i,j} x_{i',j}$ belongs to
${\rm in}_{<_{\rm rev}}(J)$.
Let $i = i'$.  
% If $k < i$, then the initial monomial of
% $\theta_{*i} x_{i,k}$ is $x_{i,k}x_{i,i}$.
Let $f$ be the polynomial 
\[
\theta_{i*}x_{i,j} - \xi_{i,j}^{-1}\theta_{*j}
\Big(\sum_{0 \leq k \leq j - 1, \, 
(i,k) \in {\mathcal A}}\xi_{i,k}x_{i,k} \Big),
\]
and write $f = f_{1} + f_{2}$, where $f_{2} \in 
{\rm in}_{<_{\rm rev}}(I_{{\mathcal O}(P)})$ and where none of the monomials
appearing in $f_{1}$ belongs to ${\rm in}_{<_{\rm rev}}(I_{{\mathcal O}(P)})$.
Then $f_{1} \in J$ and the initial monomial of $f_{1}$ is
$x_{i,j}^{2}$.  Hence $x_{i,j}^{2} \in {\rm in}_{<_{\rm rev}}(J)$.
Similarly, 
$x_{i,j} x_{i,j'} \in {\rm in}_{<_{\rm rev}}(J)$.
\, \, \, \QED

\begin{lemma}
\label{hibi:new}
For each $0 \leq i \leq p$, we write $j^{\sharp}_{i}$ for the smallest integer
for which $(i, j^{\sharp}_{i}) \in {\mathcal A}$.
For each $0 \leq j \leq q$, we write $i^{\flat}_{j}$ for the smallest integer
for which $(i^{\flat}_{j}, j) \in {\mathcal A}$.
Then $x_{i, j^{\sharp}_{i}}$ and $x_{i^{\flat}_{j},j}$ 
belong to ${\rm in}_{<_{\rm rev}}(J)$.
\end{lemma}

{\it Proof}\/.
Since 
$\theta_{i*}$ and $\theta_{*j}$ belong to $J$, their initial monomials
$x_{i, j^{\sharp}_{i}}$ and $x_{i^{\flat}_{j},j}$ 
belong to ${\rm in}_{<_{\rm rev}}(J)$.
\, \, \, \, \, \, \, \, \, \, 
\, \, \, \, \, \, \, \, \, \, 
\, \, \, \, \, \, \, \, \, \, 
\, \, \, \, \, \, 
\QED

Let ${\mathcal S}$ denote the set of square-free  monomials of $K[{\bf x}]$ of the form
\begin{eqnarray}
\label{hibi:standard}
x_{i_{1},j_{1}} x_{i_{2},j_{2}} \cdots x_{i_{r},j_{r}},
\end{eqnarray}
with each $(i_{k}, j_{k}) \in {\mathcal A} \setminus
(\{ \, x_{i, j^{\sharp}_{i}} \, : \, 0 \leq i \leq p \, \}
\cup \{ \, x_{i^{\flat}_{j},j} \, : \, 0 \leq j \leq q \, \})$
such that 
\[ 
0 < i_{1} < i_{2} < \cdots < i_{r} \leq p, \, \, \, 0 < j_{1} < j_{2} < \cdots < j_{r} \leq q,
\, \, \, r = 0, 1, 2, \ldots.
\]

\begin{figure}[tbp]
\setlength{\unitlength}{0.4mm}
\begin{center}
\begin{picture}(180,180)
%\multiput(75,0)(15,15){8}{\line(-1,1){75}}
%\multiput(75,0)(-15,15){6}{\line(1,1){105}}
%\multiput(75,0)(15,15){8}{\multiput(0,0)(-15,15){6}{\circle{3}}}
%%Prog: fig3.c
\put(75,0){\circle{3}}
\put(90,15){\circle{3}}
\put(105,30){\circle{3}}

\put(60,15){\circle{3}}
\put(75,30){\circle{3}}
\put(90,45){\circle{3}}
\put(105,60){\circle{3}}
\put(120,75){\circle{3}}
\put(135,90){\circle{3}}

\put(45,30){\circle{3}}
\put(60,45){\circle{3}}
\put(75,60){\circle{3}}
\put(90,75){\circle{3}}
\put(105,90){\circle{3}}
\put(120,105){\circle{3}}

\put(60,75){\circle{3}}
\put(75,90){\circle{3}}
\put(90,105){\circle{3}}
\put(105,120){\circle{3}}
\put(120,135){\circle{3}}
\put(135,150){\circle{3}}

\put(75,120){\circle{3}}
\put(90,135){\circle{3}}
\put(105,150){\circle{3}}
\put(120,165){\circle{3}}

\put(75,150){\circle{3}}
\put(90,165){\circle{3}}
\put(105,180){\circle{3}}

\put(75,0){\line(1,1){15}} \put(75,0){\line(-1,1){15}} \put(75,30){\line(1,-1){15}} \put(75,30){\line(-1,-1){15}} \put(90,15){\line(1,1){15}} \put(90,15){\line(-1,1){15}} \put(90,45){\line(1,-1){15}} \put(90,45){\line(-1,-1){15}} \put(105,30){\line(1,1){15}} \put(105,30){\line(-1,1){15}} \put(105,60){\line(1,-1){15}} \put(105,60){\line(-1,-1){15}} 
\put(60,15){\line(1,1){15}} \put(60,15){\line(-1,1){15}} \put(60,45){\line(1,-1){15}} \put(60,45){\line(-1,-1){15}} \put(75,30){\line(1,1){15}} \put(75,30){\line(-1,1){15}} \put(75,60){\line(1,-1){15}} \put(75,60){\line(-1,-1){15}} \put(90,45){\line(1,1){15}} \put(90,45){\line(-1,1){15}} \put(90,75){\line(1,-1){15}} \put(90,75){\line(-1,-1){15}} \put(105,60){\line(1,1){15}} \put(105,60){\line(-1,1){15}} \put(105,90){\line(1,-1){15}} \put(105,90){\line(-1,-1){15}} \put(120,75){\line(1,1){15}} \put(120,75){\line(-1,1){15}} \put(120,105){\line(1,-1){15}} \put(120,105){\line(-1,-1){15}} \put(135,90){\line(1,1){15}} \put(135,90){\line(-1,1){15}} \put(135,120){\line(1,-1){15}} \put(135,120){\line(-1,-1){15}} 
\put(75,60){\line(1,1){15}} \put(75,60){\line(-1,1){15}} \put(75,90){\line(1,-1){15}} \put(75,90){\line(-1,-1){15}} \put(90,75){\line(1,1){15}} \put(90,75){\line(-1,1){15}} \put(90,105){\line(1,-1){15}} \put(90,105){\line(-1,-1){15}} \put(105,90){\line(1,1){15}} \put(105,90){\line(-1,1){15}} \put(105,120){\line(1,-1){15}} \put(105,120){\line(-1,-1){15}} \put(120,105){\line(1,1){15}} \put(120,105){\line(-1,1){15}} \put(120,135){\line(1,-1){15}} \put(120,135){\line(-1,-1){15}} 
\put(90,105){\line(1,1){15}} \put(90,105){\line(-1,1){15}} \put(90,135){\line(1,-1){15}} \put(90,135){\line(-1,-1){15}} \put(105,120){\line(1,1){15}} \put(105,120){\line(-1,1){15}} \put(105,150){\line(1,-1){15}} \put(105,150){\line(-1,-1){15}} \put(120,135){\line(1,1){15}} \put(120,135){\line(-1,1){15}} \put(120,165){\line(1,-1){15}} \put(120,165){\line(-1,-1){15}} 
\put(90,135){\line(1,1){15}} \put(90,135){\line(-1,1){15}} \put(90,165){\line(1,-1){15}} \put(90,165){\line(-1,-1){15}} \put(105,150){\line(1,1){15}} \put(105,150){\line(-1,1){15}} \put(105,180){\line(1,-1){15}} \put(105,180){\line(-1,-1){15}}

%%end Prog
\put(90,45){\circle{6}}  \put(95,43){$\alpha_{1,2}$}
\put(75,60){\circle*{6}} \put(80,58){$x_{2,2}$}
\put(105,90){\circle{6}}  \put(110,88){$\alpha_{2,4}$}
\put(90,105){\circle*{6}} \put(95,103){$x_{3,4}$}
\put(105,150){\circle{6}} \put(110,148){$\alpha_{4,6}$}
\put(90,165){\circle*{6}} \put(95,163){$x_{5,6}$}
\thicklines
\multiput(74.7,-0.3)(0.3,0.3){3}{\line(-1,1){15}}
\multiput(60.3,14.7)(-0.3,0.3){3}{\line(1,1){30}}
\multiput(89.7,44.7)(0.3,0.3){3}{\line(-1,1){15}}
\multiput(75.3,59.7)(-0.3,0.3){3}{\line(1,1){30}}
\multiput(104.7,89.7)(0.3,0.3){3}{\line(-1,1){30}}
\multiput(75.3,119.7)(-0.3,0.3){3}{\line(1,1){30}}
\multiput(104.7,149.7)(0.3,0.3){3}{\line(-1,1){15}}
\multiput(90.3,164.7)(-0.3,0.3){3}{\line(1,1){15}}
\end{picture}
\end{center}
\caption{}  \label{fig:s57}
\end{figure}

\begin{theorem}
\label{hibi:twochainstandard}
The set of standard monomials of ${\rm in}_{<_{\rm rev}}(J)$ is equal to ${\mathcal S}$.
\end{theorem}
 
{\it Proof}\/. 
% Lemmata \ref{hibi:old} and \ref{hibi:new}
% guarantee that each standard monomial must belong to ${\mathcal S}$.  
In \cite{BJO}, it was proven that the number of standard monomials of degree $r$
coincides with the number of maximal chains of ${\mathcal J}(P)$ with $r$ descents.
Recall that the descents of a maximal chain
\begin{eqnarray*}
\label{hibi:descent}
\alpha_{0,0} = \alpha_{i_{0},j_{0}} < \alpha_{i_{1},j_{1}} 
< \cdots < \alpha_{i_{p+q},j_{p+q}} = \alpha_{p,q}
\end{eqnarray*}
of ${\mathcal J}(P)$ are those $\alpha_{i_{k},j_{k}}$ with
$1 \leq k < p+q$ such that
\[
i_{k-1} = i_{k} < i_{k+1}, \, \, \, \, \, j_{k-1} < j_{k} = j_{k+1},
\, \, \, \, \, j_{k+1} \neq j^{\sharp}_{i_{k+1}}.
\]
Now, given a square-free  monomial (\ref{hibi:standard}) of degree $r$, 
we can associate a unique maximal chain whose descents are
\[
\alpha_{i_{1}-1,j_{1}}, \alpha_{i_{2}-1,j_{2}}, \cdots, \alpha_{i_{r}-1,j_{r}},
\]
in the obvious way (see Figure \ref{fig:s57}.)

Hence the number of square-free  monomials (\ref{hibi:standard}) of degree $r$
is less than or equal to that of standard monomials of degree $r$.
On the other hand, since Lemmata \ref{hibi:old} and \ref{hibi:new} 
guarantee that each standard monomial must belong to ${\mathcal S}$,
it follows that ${\mathcal S}$ is
the set of standard monomials of ${\rm in}_{<_{\rm rev}}(J)$,
as desired.
\QED

\begin{remark} \rm
\label{hibi:parameter}
When the variables $\eta_{i*}$ and $\eta_{*j}$ are ignored
in the rational function field, we can work with
a system of parameters consisting of homogeneous elements and both
Lemma \ref{hibi:new} and Lemma \ref{hibi:old} are valid 
without modification.  
This observation is crucial to our argument of counting 
the number of standard monomials 
in the proof of Theorem \ref{hibi:twochainstandard}.
We also note that $\xi_{ij}$ may be specialized to any nonzero number
for $K={\bf C}$ without changing the claims of this section.
\end{remark}

\begin{figure}[tbp]
\begin{minipage}{0.49\hsize}
\setlength{\unitlength}{0.7mm}
\begin{center}
\begin{picture}(17,17)(-7,-2)
\put(0,0){\circle{2}}
\put(0,10){\circle{2}}
\put(10,0){\circle{2}}
\put(10,10){\circle{2}}
\put(-7,0){$a_1$}
\put(-7,10){$a_2$}
\put(13,0){$b_1$}
\put(13,10){$b_2$}
\put(0,0){\line(0,1){10}}
\put(0,10){\line(1,-1){10}}
\put(10,0){\line(0,1){10}}
\end{picture}
\end{center}
\caption{$P$}  \label{fig:4m}
\end{minipage}
%%%%
\begin{minipage}{0.49\hsize}
\setlength{\unitlength}{0.3mm}
\begin{center}
\begin{picture}(110,65)(40,-5)
\put(78,0){$00$}
\put(93,15){$01$}
\put(108,30){$02$}

\put(63,15){$10$}
\put(78,30){$11$}
\put(93,45){$12$}

\put(63,45){$21$}
\put(78,60){$22$}

%%Prog: fig4.c
\put(75,0){\circle{3}}
\put(90,15){\circle{3}}
\put(105,30){\circle{3}}

\put(60,15){\circle{3}}
\put(75,30){\circle{3}}
\put(90,45){\circle{3}}

\put(60,45){\circle{3}}
\put(75,60){\circle{3}}

\put(75,0){\line(1,1){15}} \put(75,0){\line(-1,1){15}} \put(75,30){\line(1,-1){15}} \put(75,30){\line(-1,-1){15}} \put(90,15){\line(1,1){15}} \put(90,15){\line(-1,1){15}} \put(90,45){\line(1,-1){15}} \put(90,45){\line(-1,-1){15}} 
\put(75,30){\line(1,1){15}} \put(75,30){\line(-1,1){15}} \put(75,60){\line(1,-1){15}} \put(75,60){\line(-1,-1){15}} 

%%end Prog
\end{picture}
\end{center}
\caption{${\mathcal J}(P)$}  \label{fig:5m}
\end{minipage}
\end{figure}

\begin{example}
\label{hibi:example1}
{\em
Let $P$ be the finite poset of Figure \ref{fig:4m},
and let $L = {\mathcal J}(P)$ be the distributive lattice shown in Figure \ref{fig:5m}.
Then the standard monomials of 
${\rm in}_{<_{\rm rev}}(J)$ are 
$1$; $x_{1,1}$; $x_{1,2}$; $x_{2,2}$; and $x_{1,1}x_{2,2}$. 
}
\end{example}

Let us turn to the discussion of $A$-hypergeometric systems.
Let $P_{p,q}$ denote the disjoint union of two chains $C_{p} : a_{1} < \cdots < a_{p}$
of length $p-1$ 
and $C_{q} : b_{1} < \cdots < b_{q}$ of length $q-1$.  
Let $\alpha_{i,j}$, where $0\leq i \leq p$ and $0\leq j \leq q$, be the poset ideal 
$\{a_{1},\ldots,a_{i}, b_{1},\ldots,b_{j}\}$.
In particular $\alpha_{0,0} = \emptyset$.
This is a special and interesting subclass of poset ideals.
We regard the vector $w_\alpha$, $\alpha \in {\mathcal J}(P_{p,q})$,
as a column vector and construct a matrix $A_{p,q}$ with these column vectors and a row vector
$(1,1, \ldots, 1)$.
For example, $A_{2,2}$ is
$$
\bordermatrix{
&00 & 01 & 02 & 10 & 11 & 12 & 20 & 21 & 22 \cr
&1&1&1& 1&1&1& 1&1&1 \cr 
&0&0&0& 1&1&1& 1&1&1 \cr 
&0&0&0& 0&0&0& 1&1&1 \cr 
&0&1&1& 0&1&1& 0&1&1 \cr 
&0&0&1& 0&0&1& 0&0&1 \cr 
}.
$$
By elementary row transformations, we transform
the matrix $A_{p,q}$ into the matrix ${\bar A}_{p,q}$ of the form
$(\ref{eq:aaa})$
with $k=p+1$, $n_1 = \cdots = n_k = q+1$, and
$a_i = 0 \in {\bf R}^q$ when $i \equiv 1 \  {\rm mod}\, q+1$,
$a_{i+1} = e_k \in {\bf R}^q$ when $i \equiv k \  {\rm mod}\, q+1$.
For example, $A_{2,2}$ can be transformed into
$$
{\bar A}_{2,2}=
\bordermatrix{
&00 & 01 & 02 & 10 & 11 & 12 & 20 & 21 & 22 \cr
&1&1&1& 0&0&0& 0&0&0 \cr 
&0&0&0& 1&1&1& 0&0&0 \cr 
&0&0&0& 0&0&0& 1&1&1 \cr 
&0&1&0& 0&1&0& 0&1&0 \cr 
&0&0&1& 0&0&1& 0&0&1 \cr 
}.
$$
We note that $A_{p,q}$ and ${\bar A}_{p,q}$ define the same $A$-hypergeometric system.

The matrix $A$, which represents a poset $P$ that can be decomposed into two chains (as considered in this section),
is obtained by removing some columns from ${\bar A}_{p,q}$.
For example, the $A$ that represents Figure \ref{fig:5m} is obtained by deleting the seventh column of the matrix
${\bar A}_{2,2}$.
Therefore, for $A$,
the sparse polynomials $f_j$ (\ref{eq:fdef}) 
that can be decomposed into two chains are linear
in $t$.
In particular, the $f_j$'s for $P_{p,q}$ are in the general linear position.
It follows from the integral representation (\ref{eq:kernel}) that
the ${\bar A}_{p,q}$-hypergeometric system agrees with the
Aomoto-Gel'fand system $E(p+1, (p+1)+(q+1))$ \cite{aomoto-kita},
because the matrix ${\bar A}_{p,q}$ defines a hyperplane arrangement $V(\prod t_i \prod f_j)$ 
in a general position.
The initial ideal ${\rm in}_{<_{\rm rev}}(I_A)$
is a square-free  monomial ideal.
In particular, it follows from Corollary \ref{cor:std} that
the standard monomials of the ideal $J$ defined in (\ref{eq:J})
provide a basis of the Pfaffian system for $H_A(s)$
when $x_{ij}$ is replaced by $\pd{ij}$.

Let ${\cal S}$ be the set of standard monomials given in
Theorem \ref{hibi:twochainstandard} for the poset ideal ${\mathcal J}(P)$.
Then, the set ${\cal S} |_{x_{ij} \rightarrow \pd{ij}}$ gives a basis
of the Pfaffian system for the $A$-hypergeometric system.
 
From our Theorems \ref{th:cbasis} and \ref{hibi:twochainstandard}
and the correspondence that we have explained above,
we have the following theorem.

\begin{theorem}  \label{th:twisted}
Let $A$ be the matrix representing a poset $P$ that can be decomposed into two chains,
and let ${\cal S}=\{u_1, \ldots, u_r\}$ be the set of standard monomials given in
Theorem \ref{hibi:twochainstandard} with $x_{ij}$ replaced by $\pd{ij}$.
Set $Q=\prod_{j=1}^k f_j(x,t)^{\alpha_j} t^\gamma$
and $Q' = Q |_{\alpha_i=\gamma_j = 1}$.
Then, the set of rational forms
\begin{equation} \label{eq:rational_form}
   \frac{u_i \bullet Q}{Q} dt_1 \cdots dt_m, \quad i=1, \ldots, r
\end{equation}
is a basis of the twisted cohomology group $H^{m}(\Omega^\bullet(*Q'),\nabla)$
when $\alpha_i, \gamma_j$ are generic complex numbers.
\end{theorem}

In the case $P=P_{p,q}$, this theorem is a different presentation of the celebrated work of K.~Aomoto, 
who gave a basis for the twisted cohomology group
associated with a hyperplane arrangement in a general position
(see, e.g., \cite[Theorem 9.6.2]{aomoto-kita}).
In a more general result,
Orlik and Terao gave bases of twisted cohomology groups associated with hyperplane arrangements
in terms of the $\beta nbc$ basis \cite[6.3]{orlik-terao}.
Our theorem gives bases for twisted cohomology groups in a very different way
for a class of hyperplane arrangements obtained by restricting the arrangements in the general position
to the $x_{ij}=0$'s.

\begin{figure}[tb]
\setlength{\unitlength}{5mm}
\begin{center}
\begin{picture}(2,7)(-1,-3)
\put(-1,0){\line(1,0){3}}
\put(0,-2){\line(0,1){6}}
\put(-1,-1){\line(1,1){3}}
\put(-1,2){\line(1,-1){4}}
\put(-1,5){\line(1,-2){4}}
\end{picture}
\end{center}
\caption{$V(t_1 t_2\prod f_j)$}  \label{fig:arrangement}
\end{figure}

\begin{example}\rm \label{ex:e36r}
The $A$-hypergeometric system associated with Figure \ref{fig:5m}
is the restriction of $E(3,6)$ to $x_{20}=0$.
Figure \ref{fig:arrangement} illustrates the arrangement that represents it.
\end{example}

\section{Rank of a Class of Order Polytopes}

We now turn to the discussion of the normalized volume of order polytopes.
It follows from \cite{rstan} that the normalized volume of the order polytope
${\mathcal O}(P)$ is equal to $e(P)$, the number of linear extensions of $P$.
Recall that an {\em antichain} of $P$ is a subset $B$ of $P$ such that
if $a$ and $b$ belong to $B$ with $a \neq b$, then $a$ and $b$ are incomparable 
in $P$.  The {\em width} of $P$ is the supremum of cardinalities of antichains of $P$.
The {\em length} of a chain $C$ is $|C| - 1$. 
The {\em rank} of $P$ is the supremum of lengths of chains of $P$.  

\begin{lemma}
\label{hibi:polynomial}
Fix positive integers $d$ and $r$.
Let $P$ be the disjoint union of $d$ chains $C_{1}, \ldots, C_{d}$, 
and assume that the length of each chain $C_{i}$ with $1 \leq i < d$ 
is at most $r - 1$.
Then there exists a polynomial $f(n)$ in $n$ of degree $r(d-1)$ such that
$e(P)$ is at most $f(n)$, where $n = |P|$.  
\end{lemma} 

{\it Proof}\/.
Let $\ell_{i}$ denote the length of $C_{i}$.
Then the number of linear extensions of $P$ is 
\[
e(P) = {n \choose \ell_{1},\ell_{2}, \ldots, \ell_{d}}
% {n \choose \ell_{1}}{n - \ell_{1}\choose \ell_{2}}{n - \ell_{1} - \ell_{2}\choose \ell_{3}}
% \cdots {n - \ell_{1} - \cdots - \ell_{d-1} \choose \ell_{d}}.
= \frac{n!}{\ell_{1}!\ell_{2}!\ldots \ell_{d}!}.
\] 
Since $\ell_{d} = n - \sum_{i=1}^{d-1}\ell_{i} \geq n - r(d-1)$, it follows that
\[
e(P) \leq \frac{n!}{\ell_{d}!} \leq \frac{n!}{(n - r(d-1))!}.
\]
Let
\[
f(n) = n(n-1)(n-2)\cdots(n - r(d-1) + 1),
\]
which is a polynomial in $n$ of degree $r(d-1)$.
Then $e(P) \leq f(n)$, as required.
\QED

\begin{theorem} \label{hibi:polynomialorder}
Fix positive integers $d$ and $r$.
Let $P$ be a finite partially ordered set, 
and suppose that there exists a chain $C$ of $P$ such that
\begin{enumerate}
\item[{\rm (i)}] 
the width of $P \setminus C$ is at most $d - 1$; 
\item[{\rm (ii)}] 
the rank of $P \setminus C$ is at most $r - 1$.  
\end{enumerate}
Then there exists a polynomial $f(n)$ in $n$ of degree $r(d-1)$ such that
$e(P)$ is at most $f(n)$, where $n = |P|$.
\end{theorem}

{\it Proof}\/.
Since the width of $P \setminus C$ is at most $d - 1$,
Dilworth's theorem \cite{DIL}
guarantees the existence of
$d - 1$ chains $C_{1}, \ldots, C_{d-1}$ of $P \setminus C$, 
where the length of each $C_{i}$ is at most $r - 1$,
such that $P \setminus C = C_{1} \cup C_{2} \cup \cdots \cup C_{d-1}$
and $C_{i} \cap C_{j} = \emptyset$ for $i \neq j$. 
Hence there exists a partially ordered set $Q$
that is the disjoint union of $d$ chains $C'_{1}, \ldots, C'_{d}$,
where 
the length of each $C'_{i}$ with $1 \leq i < d$ 
is at most $r - 1$,
such that there is an order-preserving bijection
$\varphi: Q \rightarrow P$.
Hence $e(P) \leq e(Q)$.  Thus the desired result follows from Lemma \ref{hibi:polynomial}.
\QED

% \bibitem{DIL}
% R. P. Dilworth, Robert P. (1950), 
% A decomposition theorem for partially ordered sets,
% {\em Annals of Math.} {\bf 51} (1950), 161 -- 166.

% \bibitem{rstan}
% R. P. Stanley,
% Two poset polytopes, {\em Discrete Comput. Geom.} {\bf 1} (1986), 9 -- 23.

Let us turn to the discussion of $A$-hypergeometric systems.
By Theorem \ref{hibi:polynomialorder}, we can regard
the rank $e(P)$ of the hypergeometric system associated with
the order polytope ${\cal O}(P)$ ($n=|P|$) has a polynomial growth property
with respect to $n$.  
This is good news, since the rank determines the complexity of the holonomic gradient method \cite{n3ost2}.

\section{Bouquet}
We now wish to introduce a ``bouquet'' of 
finite distributive lattices.
Let $P_{1}, \ldots, P_{q}$ be finite posets, where
\[
P_{i} = \{ a_{1}^{(i)}, \ldots, a_{d_{i}}^{(i)} \}, \, \, \, \, \, 
\, \, \, \, \, 
1 \leq i \leq q, 
\]
and let $L_{i} = {\mathcal J}(P_{i})$ be the distributive lattice consisting of
all poset ideals of $P_{i}$.
The finite meet-semilattice (\cite[p.~249, 3.3]{ECvol1SecEd}) 
$\bigcup_{i=1}^{q}({\mathcal J}(P_{i}))$
is called the {\em bouquet} of 
$L_{1} = {\mathcal J}(P_{1}), \ldots, L_{q} = {\mathcal J}(P_{q})$. 
For example, if $q = 3$ and $P_{1} = P_{2} = P_{3}$ are the finite poset $P_{1,1}$
shown in Figure \ref{fig:p11}, 
then the Hasse diagram of the bouquet of $L_{1}, L_{2}, L_{3}$
is shown in Figure \ref{fig:F_A}.

\begin{figure}[tbp]
\begin{minipage}{0.49\hsize}
\setlength{\unitlength}{1mm}
\begin{center}
\begin{picture}(15,23)(-15,-3)
\put(0,0){\circle{2}}
\put(10,10){\circle{2}}
\put(-10,10){\circle{2}}
\put(0,20){\circle{2}}
\put(0,0){\line(1,1){10}}
\put(10,10){\line(-1,1){10}}
\put(0,0){\line(-1,1){10}}
\put(-10,10){\line(1,1){10}}
\end{picture}
\end{center}
\caption{${\mathcal J}(P_{1,1})$}  \label{fig:p11}
\end{minipage}
%%%%
\begin{minipage}{0.49\hsize}
\setlength\unitlength{0.005cm}
\begin{picture}(650,650)(-650,-50)
\put(0,0){\circle*{20}}
\put(-10,20){$00$}
\put(0,0){\line(1,2){50}}
\put(0,0){\line(-1,2){50}}
\put(50,100){\circle*{20}}
\put(55,110){$01$}
\put(-50,100){\circle*{20}}
\put(-55,110){$10$}
\put(50,100){\line(-1,2){50}}
\put(-50,100){\line(1,2){50}}
\put(0,200){\circle*{20}}
\put(-10,220){$11$}

\put(0,0){\line(1,1){200}}
\put(0,0){\line(-1,1){200}}
\put(200,200){\line(-1,1){200}}
\put(-200,200){\line(1,1){200}}
\put(200,200){\circle*{20}}
\put(-200,200){\circle*{20}}
\put(210,210){$02$}
\put(-210,210){$20$}
\put(0,400){\circle*{20}}
\put(-10,420){$22$}

\put(0,0){\line(2,1){600}}
\put(0,0){\line(-2,1){600}}
\put(600,300){\line(-2,1){600}}
\put(-600,300){\line(2,1){600}}

\put(600,300){\circle*{20}}
\put(-600,300){\circle*{20}}

\put(0,600){\circle*{20}}
\put(610,310){$03$}
\put(-610,310){$30$}
\put(-10,620){$33$}

\end{picture}
\caption{The bouquet of $3$ ${\mathcal J}(P_{1,1})$'s}
\label{fig:F_A}
\end{minipage}
\end{figure}

Let ${\bf e}^{(i)}_{j}$, $1 \leq i \leq q, \, 1 \leq j \leq d_{i}$, denote
the standard unit coordinate vectors of $\RR^{d}$, where $d = d_{1} + \cdots + d_{q}$.
If $\beta$ is a subset of $P_{i}$, then we write $w_{\beta}$ for the $(0, 1)$-vector
$\sum_{a_{j}^{(i)}\in \beta} {\bf e}^{(i)}_{j} \in \RR^{d}$.
In particular, 
$w_{\emptyset}$ is the origin of $\RR^{d}$.
Let ${\mathcal O}(P_{1}, \ldots, P_{q}) \subset \RR^{d}$
denote the convex hull of the finite set
\[
\{ \, w_{\alpha} \, : \, \alpha \in \bigcup_{i=1}^{d}{\mathcal J}(P_{i}) \, \}.
\]
Its dimension is $\dim {\mathcal O}(P_{1}, \ldots, P_{q}) = d$.
In the language of combinatorics, the convex polytope 
${\mathcal O}(P_{1}, \ldots, P_{q})$
is called the {\em free sum} (\cite{BJM}) of
${\mathcal O}(P_{1}), \ldots, {\mathcal O}(P_{q})$.

Let 
\[
K = {\bf C}( \, \{ \, \xi_{\alpha} \, : \, \alpha \in \bigcup_{i=1}^{d}{\mathcal J}(P_{i}) \, \}, \, 
\{ \, \eta^{(i)}_{j} \, : \, 1 \leq i \leq q, \, 1 \leq j \leq d_{i} \, \}, \, \eta_{0} \, )
\]
denote the rational function field 
in $|\bigcup_{i=1}^{d}{\mathcal J}(P_{i}) + ( d + 1 )|$ variables over ${\bf C}$, and let
\[
A = K[ \, \{ \, t^{(i)}_{j} \, : \, 1 \leq i \leq q, \, 1 \leq j \leq d_{i} \, \}, \, s \, ]
\]
be the polynomial ring in $d + 1$ variables over $K$.
If $\beta$ is a subset of $P_{i}$, then we write $u_{\beta}$ for the square-free  monomial
$(\prod_{a_{j}^{(i)}\in \beta} t^{(i)}_{j})s$.
The toric ring $K[{\mathcal O}(P_{1}, \ldots, P_{q})]$ 
of ${\mathcal O}(P_{1}, \ldots, P_{q})$
is the subalgebra of $A$ that is generated by those square-free  monomials
$u_{\alpha}$ with $\alpha \in \bigcup_{i=1}^{d}{\mathcal J}(P_{i})$.
Its Krull dimension is $d + 1$.
 
Let $K[ \, \{ \, x_{\alpha} \}_{\alpha} \, ] = 
K[ \, \{ \, x_{\alpha} \, : \, \alpha \in \bigcup_{i=1}^{d}{\mathcal J}(P_{i}) \, \} \, ]$
denote the polynomial ring in $|\bigcup_{i=1}^{d}{\mathcal J}(P_{i})|$ 
variables over $K$, and define the surjective ring homomorphism 
\[
\pi : K[ \, \{ \, x_{\alpha} \}_{\alpha} \, ]
% K[ \, \{ \, x_{\alpha} \, : \, \alpha \in \bigcup_{i=1}^{d}{\mathcal J}(P_{i}) \, \} \, ] 
\to K[{\mathcal O}(P_{1}, \ldots, P_{q})]
\] 
by setting $\pi(x_{\alpha}) = u_{\alpha}$.
Its kernel is the toric ideal $I_{{\mathcal O}(P_{1}, \ldots, P_{q})}$ of 
${\mathcal O}(P_{1}, \ldots, P_{q})$.  It follows that 
$I_{{\mathcal O}(P_{1}, \ldots, P_{q})}$ is generated by those quadratic binomials
\begin{eqnarray}
\label{hibirelationagain}
x_{\alpha}x_{\beta} - x_{\alpha \wedge \beta}x_{\alpha \vee \beta},
\end{eqnarray}
where both $\alpha$ and $\beta$ belong to ${\mathcal J}(P_{i})$ 
for some $1 \leq i \leq q$ and where $\alpha$ and $\beta$
are incomparable in ${\mathcal J}(P_{i})$. 

We fix an ordering $<$ of the variables of
$K[ \, \{ \, x_{\alpha} \}_{\alpha} \, ]$ 
% $K[ \, \{ \, x_{\alpha} \, : \, \alpha \in \bigcup_{i=1}^{d}{\mathcal J}(P_{i}) \, \} \, ]$
with the property that if both $\alpha$ and $\beta$ belong to ${\mathcal J}(P_{i})$
for some $1 \leq i \leq q$
and if $\alpha > \beta$ in ${\mathcal J}(P_{i})$, then $x_{\alpha} < x_{\beta}$.
Let $<_{\rm rev}$ denote the reverse lexicographic order on
$K[ \, \{ \, x_{\alpha} \}_{\alpha} \, ]$
% $K[ \, \{ \, x_{\alpha} \, : \, \alpha \in \bigcup_{i=1}^{d}{\mathcal J}(P_{i}) \, \} \, ]$
induced by the ordering $<$.  It then follows that the set of binomials 
$(\ref{hibirelationagain})$ is the reduced Gr\"obner basis 
of $I_{{\mathcal O}(P_{1}, \ldots, P_{q})})$
with respect to $<_{\rm rev}$.
Thus the initial ideal ${\rm in}_{<_{\rm rev}}(I_{{\mathcal O}(P_{1}, \ldots, P_{q})})$ 
of $I_{{\mathcal O}(P_{1}, \ldots, P_{q})}$ with respect to $<_{\rm rev}$ is generated by 
those square-free  quadratic monomials $x_{\alpha}x_{\beta}$ such that 
both $\alpha$ and $\beta$ belong to ${\mathcal J}(P_{i})$ for some $1 \leq i \leq q$
and that $\alpha$ and $\beta$ are incomparable in ${\mathcal J}(P_{i})$.

Let
\[
\theta^{(i)}_{j} = 
\sum_{a^{(i)}_{j}\in \alpha} \xi_{\alpha}x_{\alpha} - \eta^{(i)}_{j}, \, \, \, \, \, \, \, \, \, \, 
1 \leq i \leq q, \, \, \, 1 \leq j \leq d_{i},
\]
and let
\[
\theta_{0} = \sum_{\alpha \in {\mathcal J}(P)} \xi_{\alpha}x_{\alpha} - \eta_{0}.
\]
It then follows that the sequence 
\[
(\theta_{0}, \theta^{(1)}_{1}, \dots, \theta^{(1)}_{d_{1}}, \theta^{(2)}_{1},
\ldots, \theta^{(2)}_{d_{2}}, \ldots, \theta^{(q)}_{1}, \ldots, 
\theta^{(q)}_{d_{q}})
\]
is a system of parameters of both the residue rings 
$K[\{x_{\alpha}\}_{\alpha}]/I_{{\mathcal O}(P_{1}, \ldots, P_{q})}$ 
and 
$K[\{x_{\alpha}\}_{\alpha}]/{\rm in}_{<_{\rm rev}}(
I_{{\mathcal O}(P_{1}, \ldots, P_{q})}
)$.

Now, suppose that each poset $P_{i}$ can be decomposed into two chains,
and write ${\mathcal S}_{i}$ for the set of standard monomials,
which is obtained using Theorem \ref{hibi:twochainstandard},
for the residue class ring arising from ${\mathcal O}(P_{i})$.
By virtue of the fact that no role of $\theta_{0}$ 
is required in the proof of Lemma \ref{hibi:old},
it follows that the set of standard monomials of
\begin{eqnarray}
\label{hibi:bouquetresidue}
K[\{x_{\alpha}\}_{\alpha}]/{\rm in}_{<_{\rm rev}}
({\rm in}_{<_{\rm rev}}(I_{{\mathcal O}(P_{1}, \ldots, P_{q})}),
\theta_{0}, \theta^{(1)}_{1}, 
\ldots, 
% \theta^{(1)}_{d_{1}}, 
% \ldots, \theta^{(q)}_{1}, \ldots, 
\theta^{(q)}_{d_{q}})
\end{eqnarray}
with respect to $<_{\rm rev}$ is a subset of 
\begin{eqnarray}
\label{hibi:bouquetstandard}
\Big\{ \, \prod_{i=1}^{q} u_{i} \, : \, u_{i} \in {\mathcal S}_{i}, \, 1 \leq i \leq q \, \Big\}.
\end{eqnarray}

Finally, the computation of the number of standard monomials
based on the equality $(7)$ of \cite[Theorem 1.4]{BJM} together with
the information on the facets of the order polytopes 
(\cite[p.~10]{rstan})
guarantee the following theorem.

\begin{theorem}
\label{hibi:bouquettheorem}
The set of standard monomials of the residue class ring
(\ref{hibi:bouquetresidue}) with respect to $<_{\rm rev}$
coincides with the set of square-free  monomials
(\ref{hibi:bouquetstandard}).
\end{theorem}

\begin{example} \rm
For the bouquet of Figure \ref{fig:F_A}, the set of standard monomials obtained 
by Theorem \ref{hibi:bouquettheorem} (labeling variables as in Figure \ref{fig:F_A} and replacing $x_{ij}$ by $\pd{ij}$) is
$$
\{ \pd{11}^{k_1} \pd{22}^{k_2} \pd{33}^{k_3} \,|\,
    k_i \in \{ 0, 1 \} \}.
$$
We will show that this bouquet represents the Lauricella 
hypergeometric function $F_A$ of three variables \cite[Chapitre VII]{ak}.
We note that the twisted cohomology groups associated with the $F_A$ are studied in \cite{matsu}
in a quite different way.
We consider $A$-hypergeometric system associated with 
the lattice shown in Figure \ref{fig:F_A}.
%%Prog: A-hg/Prog/toric2.rr, ttest2()
The independent variables of the system will be denoted by
$ p_{00}, p_{01}, p_{02}, p_{03}, p_{10}, p_{20}, p_{30}, p_{11}, p_{22}, p_{33}$
or simply as
$ 00, 01, 02, 03, 10, 20, 30, 11, 22, 33$ if no confusion arises.
The differential operators of the system will be denoted by
$ \pd{00}, \pd{01}, \pd{02}, \pd{03}, \pd{10}, \pd{20}, \pd{30}, \pd{11}, \pd{22}, \pd{33}$
or simply as	
$ 00, 01, 02, 03, 10, 20, 30, 11, 22, 33$ if no confusion arises.
The toric ideal associated with the lattice is generated by
\begin{equation}
\underline{10 \cdot 01} - 00 \cdot 11,
\underline{20 \cdot 02} - 00 \cdot 22,
\underline{30 \cdot 03} - 00 \cdot 33. 
\end{equation}
The underlined terms are the leading terms for the reverse lexicographic order
such that $00 < \mbox{other variables}$,
and the set is a Gr\"obner basis with this order.
Hence,
the $A$-hypergeometric system has a solution of the form
\begin{equation}
p^\gamma f\left( \frac{10 \cdot 01}{00 \cdot 11},
                      \frac{20 \cdot 02}{00 \cdot 22},
                      \frac{30 \cdot 03}{00 \cdot 33}
\right).
\end{equation}
Set
\begin{eqnarray*}
 x &=&  \frac{10 \cdot 01}{00 \cdot 11}, \\ 
 y &=&  \frac{20 \cdot 02}{00 \cdot 22}, \\
 z &=&  \frac{30 \cdot 03}{00 \cdot 33}.
\end{eqnarray*}
The differential operator
$p_{10}p_{01}p_{00}p_{11} (\underline{\pd{10}\pd{01}} - \pd{00}\pd{11})$
can be written as
$ \theta_{10} \theta_{01} - x \theta_{00} \theta_{11}$,
where $\theta_{ij} = p_{ij} \pd{ij}$ (the Euler operator).
We will derive a differential operator that annihilates the function $f$
from this operator.
Apply $\theta_{11}$ to $p^\gamma f(x,y,z)$.
Then, we have
$p^\gamma (\gamma_{11} - \theta_x) f$.
Apply $\theta_{00}$ to this function.
Then, we have
$$ 
 \gamma_{00}\gamma_{11} p^\gamma f 
+ p^\gamma \gamma_{11} ( -\theta_x-\theta_y-\theta_z) f
- \gamma_{00} p^\gamma x f - p^\gamma(-1)x f_x
- p^\gamma x  ( -\theta_x-\theta_y-\theta_z) f_x.
$$
This can be factored as
$$ p^\gamma ( \theta_x+\theta_y+\theta_z-\gamma_{00}) 
                  (\theta_x - \gamma_{11}) f.
$$
An analogous calculation leads us to
$$ \theta_{10} \theta_{01} p^\gamma f(x,y,z) =
 p^\gamma (\theta_x+\gamma_{01}) (\theta_x+\gamma_{10}) f.
$$
Therefore, the function $f(x,y,z)$ satisfies
\begin{equation}
\left( (\theta_x+\gamma_{01}) (\theta_x+\gamma_{10})  
 - x  ( \theta_x+\theta_y+\theta_z-\gamma_{00}) (\theta_x - \gamma_{11}) 
\right) f = 0.
\end{equation}
By analogous calculations, we have
\begin{eqnarray}
&&\left( (\theta_y+\gamma_{02}) (\theta_y+\gamma_{20})  
 - y  ( \theta_x+\theta_y+\theta_z-\gamma_{00}) (\theta_y - \gamma_{22}) 
\right) f = 0, \\
&&\left( (\theta_z+\gamma_{03}) (\theta_x+\gamma_{30})  
 - z  ( \theta_x+\theta_y+\theta_z-\gamma_{00}) (\theta_z - \gamma_{33}) 
\right) f = 0. 
\end{eqnarray}
By these equations for the function $f$,
we conclude that the function
$g = x^{-\gamma_{01}} y^{-\gamma_{02}} z^{-\gamma_{03}} f(x,y,z)$
satisfies the differential equations for the Lauricella function $F_A$, $n=3$.

The bouquet of $n$ squares stands for 
the Lauricella $F_A$ of $n$ variables.
\end{example}

\bigbreak
{\it Acknowledgments}:
The authors are grateful to Prof. Kiyoshi Takeuchi
for valuable explanations of rapid decay cycles and for Remark \ref{rem:rapid_decay}, the statement that these cycles are independent with respect to the parameter $x$.

\end{document}